\documentclass[11pt,fleqn,amstex]{article}

\usepackage{latexsym}

\usepackage{amsmath}
\usepackage{amsthm}
\usepackage{amssymb}
\usepackage{amsfonts}
\usepackage{color}
\usepackage{tikz, lipsum, lmodern}
\usepackage{graphicx}
\usepackage[rightcaption]{sidecap}

\newcommand{\ba}{\begin{array}}
\newcommand{\ea}{\end{array}}

\newcommand{\be}{\begin{equation}}
\newcommand{\ee}{\end{equation}}

\newtheorem{prop}{Proposition}[section]
\newtheorem{Th}[prop]{Theorem}

\numberwithin{equation}{section}

\begin{document}

\title{\bf Forced extension of GNI techniques to dissipative systems}
\author{  Artur Kobus\thanks{E-mail: a.kobus@uwb.edu.pl}
\\ {\footnotesize Uniwersytet w Bia\l ymstoku, Wydzia{\l} Fizyki}
\\ {\footnotesize ul.\ Cio\l kowskiego 1L,  15-245 Bia\l ystok, Poland}
}

\date{ }

\maketitle

\begin{abstract}
We propose new concept of energy reservoir and effectively conserved quantity, what enables us to treat dissipative systems along the lines of the framework of Geometric Numerical Integration. Using this opportunity, we try to confirm numerically if our idea is useful. Numerical experiments show good qualitative behavior of integration technique for ODEs based on non-potential Hamiltonian formalism. It occurs that rising accuracy is a difficult task due to dissipative form of the system under scrutiny.
\end{abstract}

{\bf{Keywords:}} Geometric Numerical Integration; Discrete Gradient Method; Dissipative Systems; Hamiltonian Mechanics; Reservoir Variables.

{\bf{PACs numbers:}} 45.10.-b, 02.60.Cb, 02.70.-c, 02.70.Bf

\section{Introduction}

Classical paradigm of numerical analysis of ODEs is to find one or multiple packages that can solve well-posed problem in finite time with demanded accuracy (see e.g. \cite{I, MS, PTVF}). As opposed to this, not so long ago there occurred rising need for preserving qualitative features of ODEs exactly, when accuracy went further from our main interests. This gave birth to the paradigm of Geometric Numerical Integration (or GNI, for short, see e.g. \cite{BC, HLW}), which caused many peculiar classes of ODE integrators to pop up.

One especially interesting case of a one-step algorithm is the so-called discrete gradient scheme and preserving quantities exactly is a built-in feature of the method \cite{MQR, QC}. While discrete gradient family of methods is itself of huge interest, it is also possible to approach conservative problems in a little less direct way, using e.g. symplectic schemes that play with different, numerically induced conserved quantities while preserving the symplectic property \cite{CS, Y}.

While for conservative systems there is a plenty of disposable integrators, there is very little (if not none) algorithms designed specifically to grasp correctly dissipative  behavior. This is due to the lack of e.g. conserved quantities, although mentioned discrete gradient scheme can also recreate proper behavior of a system which energy is described by Lyapunov function \cite{MQR}. The problem is that dissipative behavior is often more complicated than this; non-conservative systems exhibit plenty of non-linear phenomena like intermittency or appearance of stable limit cycles \cite{CRS}, which are hard to describe, nonetheless extremely useful.

In this paper we try to fix the situation of dissipative systems by introduction of non-potential Hamiltonian formalism. After brief remarks concerning mainly notation we give its basic theoretical description in section $3$ and discretize it in section $4$ with checking some of its basic features. Next we concentrate our efforts on showing, by our discrete gradient procedure and some classical integrators, that our approach yields correct results for the case of damped harmonic oscillator. This happens in section $5$. Section $6$ is devoted to concluding remarks and future perspectives. 

\section{Numerical glossary}

For the sake of undisturbed comprehension we give all the indispensable definitions in convenient notation used in this paper.

We begin with brief recap of errors occurring in numerical analysis: we consider as \emph{global error} the object
\be
{\bf{e}}_{x,i} = x_i - x(t_i),
\ee
where the error is estimated for quantity $x$, in the $i^{th}$ step of numerical method, while the $i^{th}$ term in solution sequence is corresponding to  exact solution $x(t)$ in the moment $t_i$. Here we stipulate, that we use constant time-step $h$, so that $t_i=t_0+ih$.

With respect to thus obtained time-grid, we measure also the \emph{local error}
\be
\mathbb{T}_{x,i+1} = \frac{x(t_{i+1})-x(t_i)}{h} - \Phi(t_i,x(t_i);h)
\ee
with subscripts understood as previously, and $\Phi$ is the \emph{numerical flow} of the considered method (in this paper we will be concerned with one-step schemes only). It is worth stressing that the local error is method-specific.

The highest order of the term in numerical flow that agrees exactly with exact flow of considered system is called \emph{theoretical order}.

We will use local error to determine order of the method under scrutiny. Let us assume we have some scheme of  theoretical order $p$ (not to be confused with momentum! It should be clear from context), then
\be
\ba{l}
\mathbb{T}_{x,i} \approx \frac{x(t_i)+h x'(t_i)+\ldots +\frac{h^p}{p!} x^{(p)} (t_i) + O(h^{p+1})-x(t_i)}{h} - \Phi(t_i,x(t_i);h) =\\[2 ex]
= \qquad O(h^p)=c h^p
\ea
\ee
where $c$ is some constant (although it might depend on $p$). Thus determined order will be referred to as \emph{empirical order} of the method. Theoretical order is obviously a global property of the method, but empirical order is not. From now on, where the distinction needed, we will use $p_t$ to denote theoretical order and $p_e$  to denote empirical order.

From here we gain upper bound on a logarithm of local error
\be
\label{bnd}
\textrm{log max}_i |\mathbb{T}_{x,i}| \leq \textrm{log c} + p \textrm{log} h.
\ee

We use this bound to determine the order of the method by running it several times with different time-steps, and then performing linear regression on collected data. The directional constant of the straight line approximately equals $p_e$.

Some caveat is in order. We deliberately choose some base time-grid, generated with the time step $h_0$. Then we apply the numerical scheme with various time-steps, with local errors calculated for each point on the base time-grid, with the third argument of the numerical flow being current time step. In this way we obtain comparable results, on equi-grid point set.

\section{Non-potential Hamiltonian systems}

We begin with the notion of Newton's equation of motion expressed in simple second order autonomous ODE form (for the one-dimensional system)
\be
\ddot{q} = F(q),
\ee
and this equation, as usual, may be cast in the Hamiltonian form
\be
\ba{l}
\dot{q} = p,\\
\dot{p} = F(q).
\ea
\ee

Flow of these equations possesses a conserved quantity
\be
E=T+V=\frac{1}{2} p^2 -\int_{q_0}^q F(q) d q,
\ee
which, expressed exclusively in terms of coordinate $q$ and momentum $p$ (as beyond) is called Hamiltonian of the system. Despite its nice feature, that it is preserved during the time evolution of the system, it is also the generator of the equations of motion through simple differentiation, namely
\be
\ba{l}
\dot{q}=\frac{\partial H}{\partial p} =p,\\
\dot{p} = -\frac{\partial H}{\partial q}  = -V'(q)= F(q),
\ea
\ee
where the potential function is defined to be
\be
V(q) = - \int_{q_0}^q F(q) d q
\ee
and to recover the force from the potential we differentiated with respect to upper limit of the integral. {\bf{From now on we will accept this formal operation as defining the force exerted on a system through differentiation}}.

Let us introduce the dissipative force of non-potential form which we signify by $D(q,p)$. We assume that its expression is already consistent with possible constraints put on the system, hence it is given in terms of generalized coordinate and momentum. Of course, appearance of such an object would prevent from occurrence any conservative behavior, unless we proceed carefully enough, to finally include $D$ in the description of the system, so the work done by this force is considered positive.

Now let us ponder
\be
\ba{l}
\label{rom}
\dot{q}=p,\\
\dot{p}=F(q)-D(q,p),
\ea
\ee
where we begin to use a reservoir variable
\be
\label{w}
w(q,q_0,p_0)=\int_{q_0}^q D(q,p) d q= \int_{t_0}^t D(q(t),p(t)) p(t) d t
\ee
which is physically measuring the work done by dissipative forces (as the symbol $w$ suggests), and so that
\be
\label{wdot}
\dot{w} = D(q,p) p.
\ee

Thus defined quantity will be referred to as artificial integral variable, name emphasizing it does not follow usual, differential evolution.

We will benefit using the second, re-parameterized form of reservoir variable, since it does not cause any trouble with unique correspondence of $q,p$ solutions.

Now we define the non-potential Hamiltonian to be
\be
K(q,p;q_0,p_0)=\frac{1}{2} p^2 + V(q) + w(q,q_0,p_0).
\ee

Above definition gives us simple way of understanding the physical meaning of a reservoir - usually we would say that, for example, friction dissipates energy producing heat. Here we pull back this quantity into the system under scrutiny, so that it counts as a positive increment to the total energy.

The dependence of function $K$ generating equations of motion on initial conditions characterizes dissipative systems. This is a formal reflection of lack of the time-translation symmetry.

Formally it is clear, but a little bit {"}odd{"} statement, that we consider Hamiltonian with added integral term, provoked by the appearance of non-potential, dissipative force in {"}Hamilton's equations{"} (\ref{rom}). It is justified by providing a full force exerted on a system only by means of differentiation of potential together with a reservoir. We have
\be
\label{wpoq}
\frac{\partial w}{\partial q} = \frac{1}{\dot{q}} \frac{d w}{d t}= D(q,p),
\ee
restoring equations of motion in full capacity.. Differentiation with respect to $p$  is trivial, because
\be
\frac{\partial w}{\partial p}=0,
\ee
since we treat $w$ as $w(q)$.

\begin{Th}
The quantity $K$ is conserved during the time evolution of the system.
\end{Th}

\begin{proof}
\be
\dot{K}= \frac{\partial H}{\partial q} \dot{q} + \frac{\partial H}{\partial p} \dot{p} + \dot{w}= (-\dot{p}-D(q,p)) \dot{q} +\dot{q} \dot{p} +D(q,p) p =0.
\ee
 Note that the crucial part here is to exclude possibility of $D$ depending explicitly on time, this would make the problem non-autonomous.
\end{proof}

Further we will refer to introduced conserved quantity as {\bf{effectively conserved}}, the name meaning that its behavior results from equations of motion after adding a reservoir to the system, not from the equations of motion solely.

Physical interpretation of the non-potential Hamiltonian is strikingly simple: it is initial energy of the system. We can view the fact of its preservation as just kinetic energy being transformed in a two-fold way: as usual, it becomes stored in potential energy form, or it is being irreversibly {"}eaten{"} by the reservoir (it is the case only in the damped case, generic form of non-potential forces is such that it can stimulate the motion, or mutually absorb and inject energy of the reservoir $w$ into the system). 

Especially interesting is fact that we would not use new variable $w$ while solving differential equations, but it is of key importance for preserving $K$. {\bf{The main idea is simple: by considering reservoir $w$, we push back the system into an effectively conservative form}}.

As indicating from stated remarks, we lean on assumptions:
\begin{enumerate}
\item We conceive of dissipative forces as contained in the system and consistent with all the constraints, so expressed by generalized coordinate and momentum. This last variable remains uninfluenced by inclusion of additional elements in the system.

\item Dynamically, we adjoin to the system the reservoir $w$, containing work done by dissipative forces (it certainly plays no role in solution of equations of motion, thus is just a redundant variable). As an effect, the new generator of equations of motion, $K$, is conserved.

\item If we lay $D(q,p) \equiv 0$, system goes back to its pure Hamiltonian, potential form.
\end{enumerate}

\section{Modified discrete gradients}

Given equations of motion (we consider, for the time being, only one-dimensional systems - generalization to more degrees of freedeom, as notationally little cumbersome, will be handled elsewhere)
\be
\ba{l}
\dot{q} = \frac{\partial H}{\partial p},\\
\dot{p} = - \frac{\partial H}{\partial q}-D(q,p),\\
\dot{w} = D(q,p)p,
\ea
\ee
where the quantity
\be
K=H+w
\ee
is effectively conserved ($H$ is ordinary hamiltonian of conservative form), we discretize them, due to procedure of discrete gradient method \cite{QC}, but with $w$ variable changed every time $q$ changes (and it is understood that $w$ is one of arguments of $K$).

In other words, we put
\be
\label{modgr}
\ba{l}
\frac{q_{i+1}-q_i}{h} = \frac{K(q_{i+1},p_{i+1},w_{i+1}) - K(q_{i+1},p_i,w_{i+1})}{p_{i+1}-p_i},\\
\frac{p_{i+1}-p_i}{h} = \frac{K(q_i,p_{i},w_i) - K(q_{i+1},p_{i},w_{i+1})}{q_{i+1}-q_i},\\
\frac{w_{i+1}-w_i}{h} = \frac{1}{2} D(q_i,q_{i+1},p_i,p_{i+1}) (p_i + p_{i+1})
\ea
\ee
where we are able to express evolution of $w$ in quite arbitrary way - it should only obey the condition of becoming $D(q,p) p$ in the continuous case.

Above scheme guarantees that
\be
\label{pres}
K(q_i,p_i,w_i) = K(q_{i+1},p_{i+1},w_{i+1})
\ee
as declared before.

It is worth emphasis that when dissipative forces are absent, this becomes usual discrete gradient method.

In the following, we will use simple iteration technique to solve implicit equations with tolerance $\varepsilon = 10^{-18}$, base time-step will take the value $h_0=0,001$.

During the measurement we use set of time-steps: 
\be
h =\{0.036,0.03,0.02,0.028,0.017,0.01\},\\
\ee
where small range is dictated by the will to capture linear behavior during regression. Using different set would give different results in another local area. From this remark we may conjecture that empirical order is a local quantity.

Now we are ready to deal with rising the order of this gradient scheme using example of damped oscillator. We begin with equation of motion for evolution of $q$. We introduce function $\delta$ allowing us to rise the order, so that exact preservation property (\ref{pres}) would not be altered. We expand both sides of equation in power series in $h$, like in \cite{CR}, assuming that usage of $\delta_{[\infty]}$ provides access to exact integrator. We have two main options:

{\bf{1. Equation for evolution of $q$:}}

We begin with
\be
(q(t_{i+1})-q(t_i))=\frac{1}{2} \delta_{[\infty]}^{(q)}  (p(t_{i+1})+p(t_i))
\ee
which becomes
\be
\dot{q} h+ \frac{1}{2} \ddot{q} h^2 + \frac{1}{6} q^{(3)} h^3 + \ldots = \frac{1}{2} (\delta_1^{(q)} h + \delta_2^{(q)} h^2 + \ldots) (2 p + \dot{p} h + \frac{1}{2} \ddot{p} h^2 \ldots)
\ee
yielding $\delta$ coefficients
\be
\ba{l}
\delta_1^{(q)} =1,\quad \delta_2^{(q)} =0,\quad \delta_3^{(q)}= -\frac{1}{12} \frac{\ddot{p}}{p} = \frac{1-b^2-b q/p}{12},\\
\delta_4^{(q)} = \frac{-\frac{1}{24} p^{(3)}-\frac{1}{2} \dot{p} \delta_3}{p}=\frac{(q+bp)(-b^2 p - b q + p)}{24 p^2}.
\ea
\ee

In the case of conservative system we would expet normally rising the order of the gradient method with addition of every extra $\delta_i$ term. For dissipative systems, however, it turns out that theoretical order upgrade does not mean rise of empirical order (although it can). This is because we divide every coefficient by the expression that takes zero values at some points. It does not kill convergence of the method, since at the same time we multiply those coefficients by higher and higher powers of time-step. For this reason we find order of the method relatively low, compared to what is should be.

\begin{minipage}[c]{0.495\textwidth}

\includegraphics[width=6cm, height=5cm]{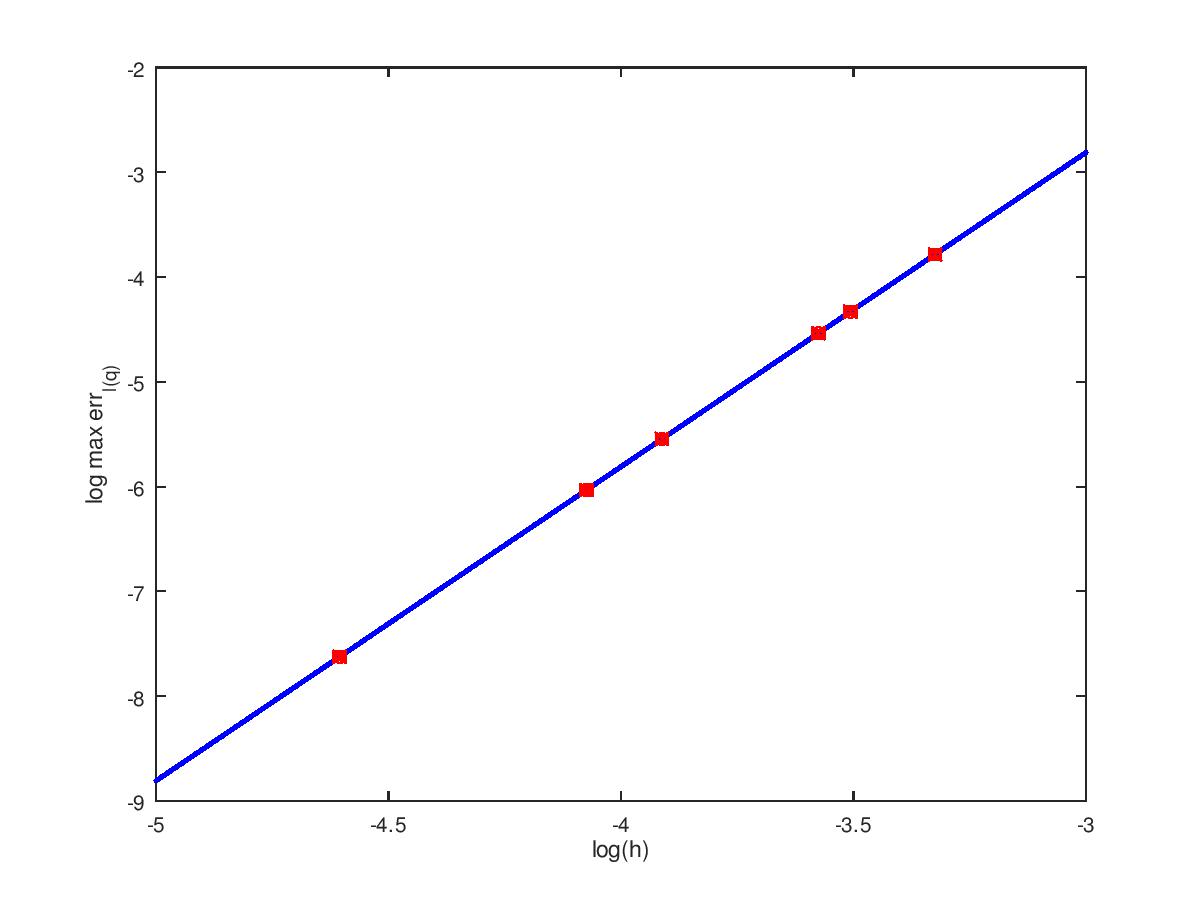}

\footnotesize{Figure 3.1.: Order of the method evaluated by linear regression from maxima of local $q$ errors in $\delta_3^{(q)}$. We obtain $p \approx 2.99811$.}
$\\[1ex]$

\end{minipage} \begin{minipage}[c]{0.01\textwidth}   \end{minipage} \begin{minipage}[c]{0.495\textwidth}

\includegraphics[width=6cm, height=5cm]{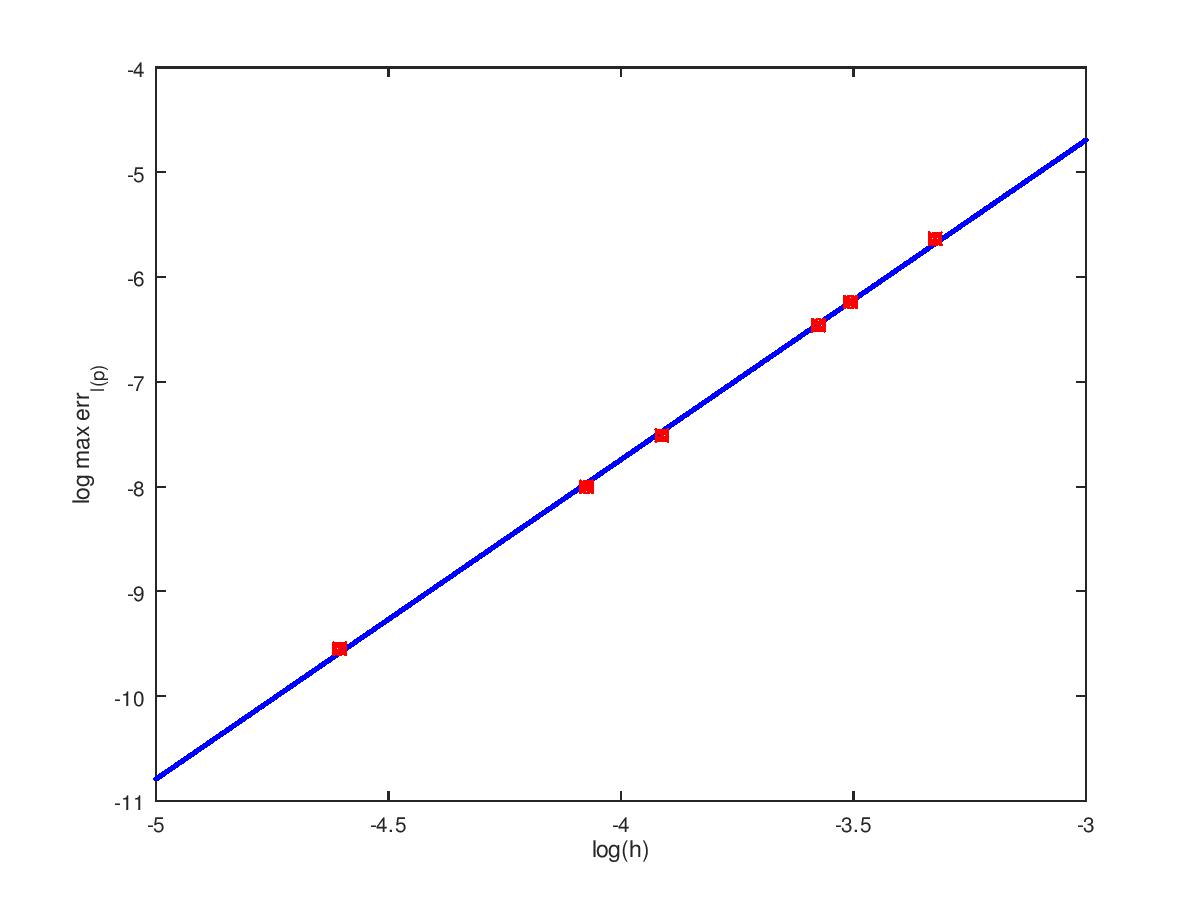}

\footnotesize{Figure 3.2.: The same numerical experiment for $\delta_3^{(q)}$, read from $p$ errors. This time order is determined to be $2.98861$.}
$\\[1ex]$

\end{minipage}

At the same time empirical order of $w$ is $2$, hence whole method is of order $2$. Addition of $\delta_4^{(q)}$ hardly improves our situation.

{\bf{2. Equation of evoltion for $p$:}}

Here we start with
\be
(p(t_{i+1})-p(t_i))=\frac{1}{2} \delta_{[\infty]}^{(p)} (- (q(t_i)+q(t_{i+1}))-b(p(t_{i+1})+p(t_i))),
\ee
where we used explicit form of terms in the equation. Hence
\be
\dot{p} h + \frac{1}{2} \ddot{p} h^2 + \ldots = -\frac{1}{2} (\delta_1^{(p)} h + \delta_2^{(p)} h^2 + \ldots) (2 q + \dot{q} h + \ldots+ b(2p + \dot{p} h + \ldots))
\ee
and we get coefficients
\be
\ba{l}
\delta_1^{(p)} =1,\quad  \delta_2^{(p)} =0,\quad \delta_3^{(p)}= \frac{1}{12} \frac{p^{(3)}}{q+bp} = \frac{-b^3 p - b^2 q + 22bp+q}{12 (q+bp)},\\
\delta_4^{(p)} = \frac{\frac{1}{24} p^{(4)}+\frac{1}{2} \ddot{p} \delta_3}{q+bp}=\frac{-b^5 p^2-2b^4qp+3b^3p^2-b^3q^2+4b^2qp+bq^2-2bp^2-qp}{24(b^2 p^2 +2bqp+ q^2)}+\\
\qquad +\frac{b^4p+b^3q-3b^2p-2bq+p}{24(bp+q)},
\ea
\ee
where we are not explicitly writing higher order coefficients since they get monstrous quickly. Note that $\delta_1=1$ guarantees consistency and $\delta_2=0$ assures we get the second order scheme at least, as previously.

Again, from the above calculation we clearly see that increasing the order of the scheme in that way should be very hard, if not impossible: with every appearance of additional power of the time-step, there occurs also additional division by $q+bp$ which at some points will cause the coefficients to blow-up. In this way we should obtain the scheme with safe second order behavior, but not higher.

Numerical experiment shows that indeed, local error committed by the method is at the second-order level, but when we run the whole procedure of determining order, we see it is growing as expected!

\begin{minipage}[c]{0.495\textwidth}

\includegraphics[width=6cm, height=5cm]{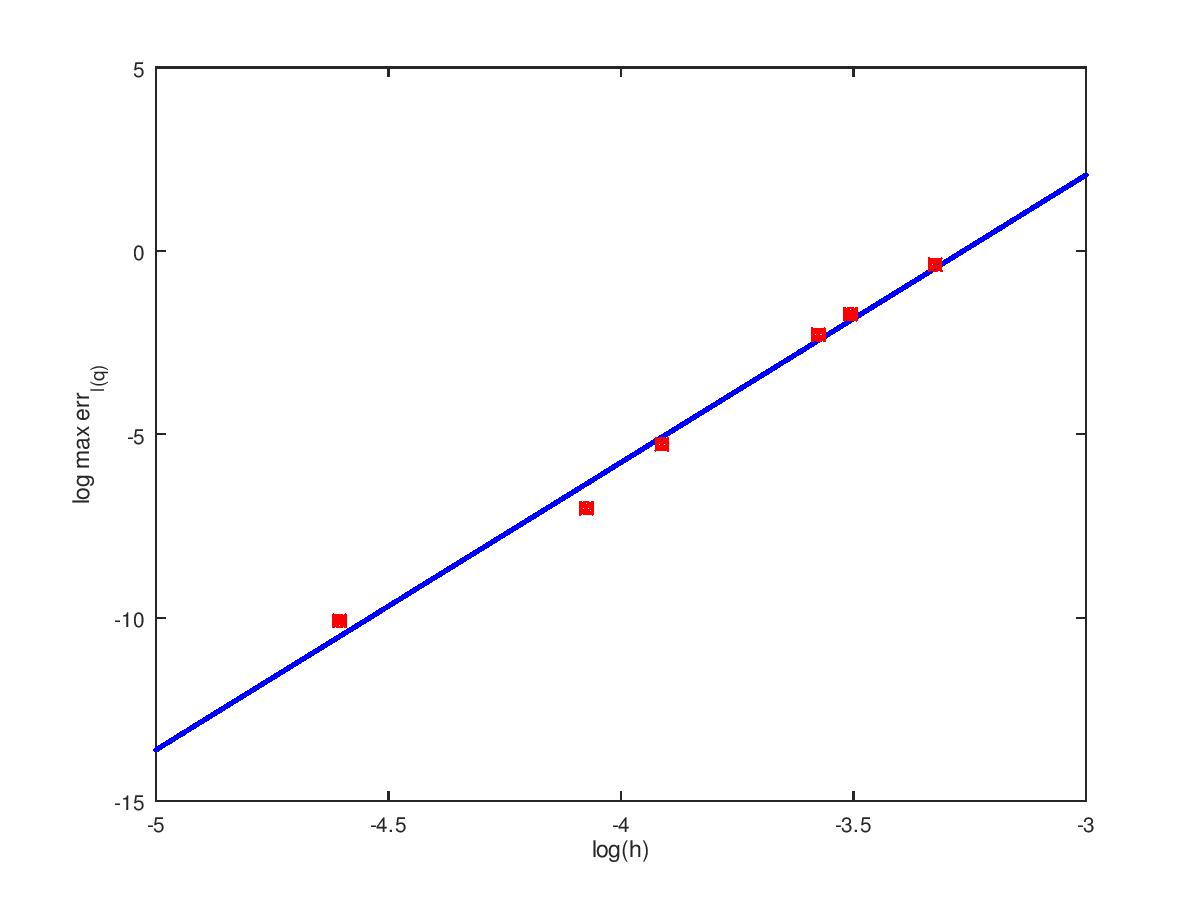}

\footnotesize{Figure 3.3.: Order of the method evaluated by linear regression from maxima of local $q$ errors in $\delta_4^{
(p)}$. We obtain $p \approx 7.83681$.}
$\\[1ex]$

\end{minipage} \begin{minipage}[c]{0.01\textwidth}   \end{minipage} \begin{minipage}[c]{0.495\textwidth}

\includegraphics[width=6cm, height=5cm]{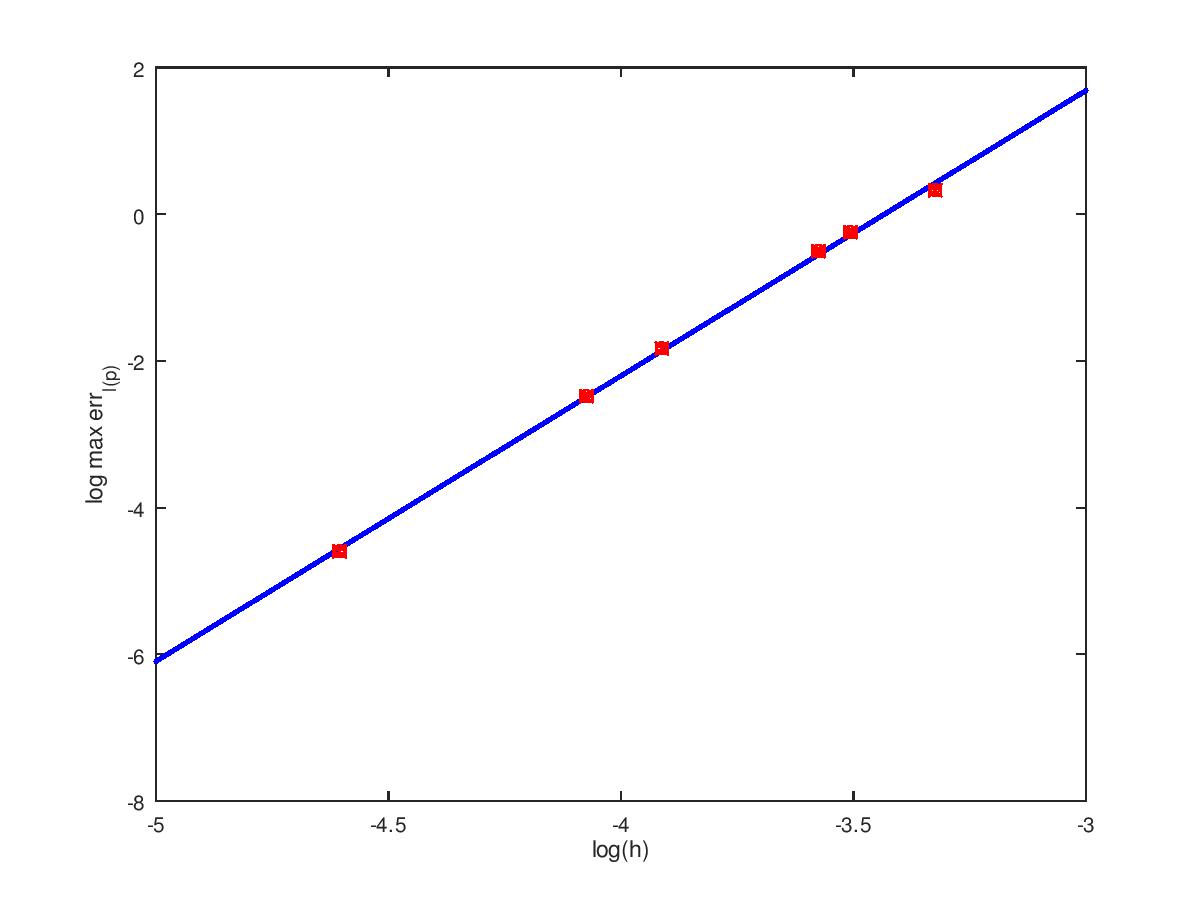}

\footnotesize{Figure 3.4.: The same numerical experiment for $\delta_4^{(q)}$, read from $p$ errors. This time order is determined to be $2.3.89048$.}
$\\[1ex]$

\end{minipage}

$w$ variable admits here behavior of order $p_e =7.7702$, so order of the method is $p_e=p_t =4$.

Changing the form of numerical evolution of $w$ would cause only a slight shift in results. Calculating $\delta$ from the third equations of continuous system meets similar problems (even appearing already in $\delta_2^{(w)}$). The reason for such behavior is exponential growth of $c$ constant with order $p_e$ in (\ref{bnd}).

\section{Numerical schemes argument}

In order to check how our new method works, we perform numerical experiment, consisting in executing few different algorithms on the same set of initial data. We compare our modification of discrete gradient method (modDG) with symplectic  leap-frog scheme (pqpLF) and explicit fourth-order Runge-Kutta (eRK4).

We use initial conditions $q_0 = 2.3, p_0=-3.1,w_0=0.0$.

Continuous system is
\be
\ba{l}
\dot{q}=p,\\
\dot{p}=-q-bp,\\
\dot{w} = bp^2,
\ea
\ee
where we stick to the caseof $k=1$, $b=0.1$ is the damping constant and we have already included reservoir in the description.

As eRK4 and modDG are clear in use with reservoir variable, the SV scheme needs a little explanation. Instead of using normal Hamiltonian, we use the $K$ generator with described earlier differentiation rules \cite{AK}. Thus
\be
\ba{l}
p_{i+\frac{1}{2}}= p_i - \frac{h}{2} (\nabla_q K(q_i,p_{i+\frac{1}{2}})) \rightarrow p_{i+\frac{1}{2}} = \frac{p_i-\frac{h}{2} q_i}{1+\frac{hb}{2}},\\
q_{i+1}=q_i+\frac{h}{2} (\nabla_p K(q_i,p_{i+\frac{1}{2}})+\nabla_p K(q_{i+1},p_{i+\frac{1}{2}})) \rightarrow q_{i+1}=q_i+hp_{i+\frac{1}{2}},\\
p_{i+1}=p_{i+\frac{1}{2}} - \frac{h}{2} (\nabla_q K(q_{i+1}, p_{i+\frac{1}{2}})) \rightarrow p_{i+1} = p_{i+\frac{1}{2}} - \frac{h}{2} (q_{i+1} + b p_{i+\frac{1}{2}})
\ea
\ee
so it is an explicit scheme.

\begin{minipage}[c]{0.495\textwidth}

\includegraphics[width=6cm, height=5cm]{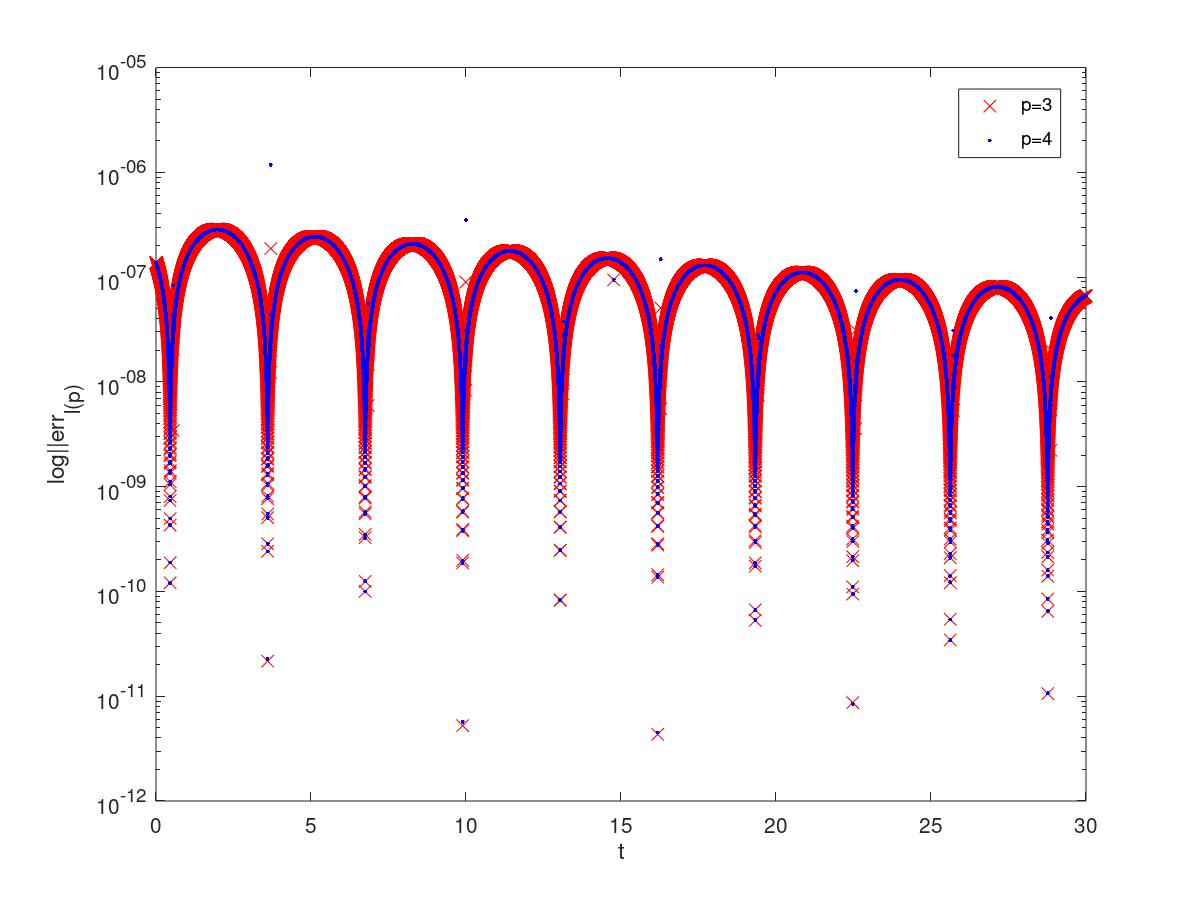}

\footnotesize{Figure 4.1.: No changes in local error of $q$ due to rising order.}
$\\[1ex]$

\end{minipage} \begin{minipage}[c]{0.01\textwidth}   \end{minipage} \begin{minipage}[c]{0.495\textwidth}

$\\[1ex]$

Picture on the left shows local errors of $q$ variable compared in two cases: when $p=4$ and we keep four terms in $\delta$ function (blue points), and when $p=3$ and we keep three terms in $\delta$ function (red crosses) so we see, that errors remain of constant magnitude. Behavior of global errors yields the same pattern. In the same time order practically does change.

$\\[1ex]$

\end{minipage}

Of course, when $p$ stays fixed smaller $h$ means smaller local error. 

\begin{minipage}[c]{0.495\textwidth}

\includegraphics[width=6cm, height=5cm]{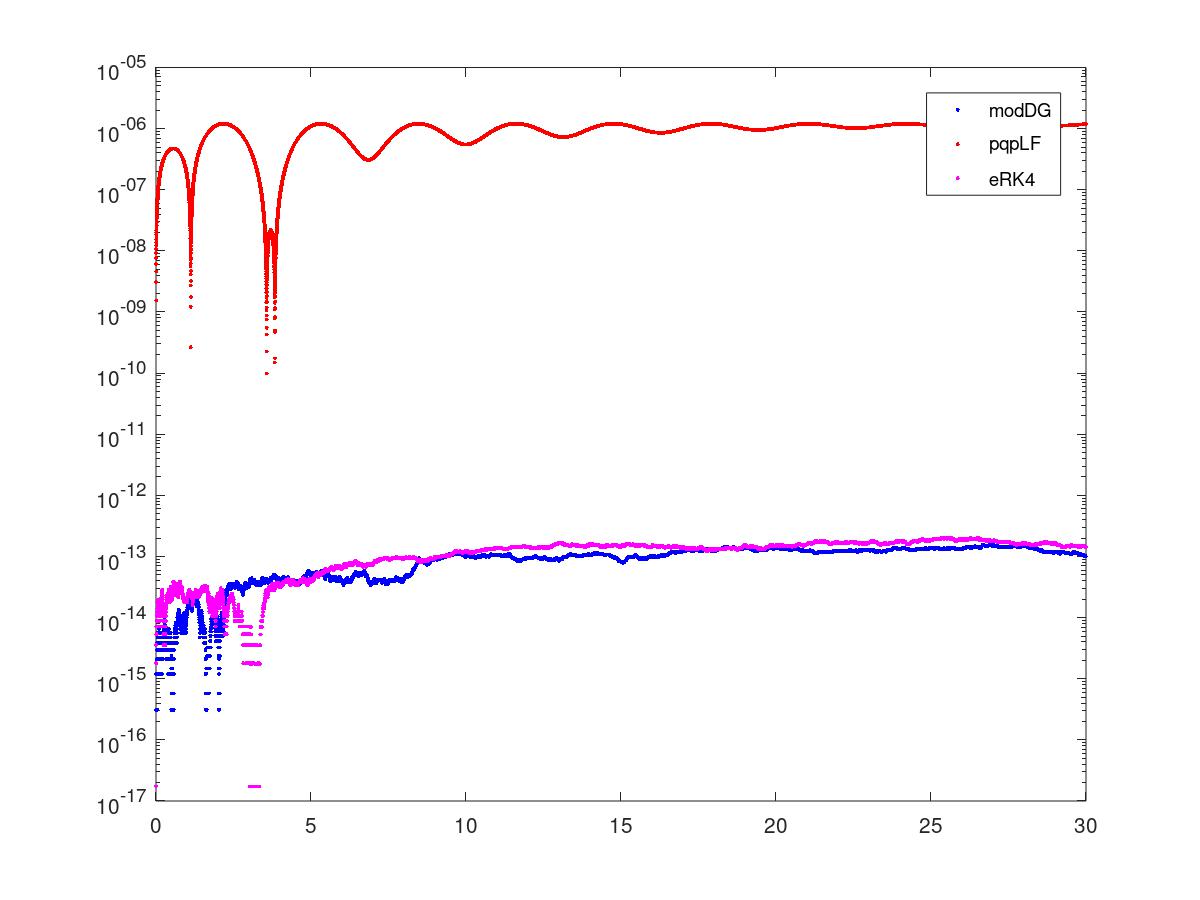}

\footnotesize{Figure 4.2.: Deviation from initial value of $K$. SV scheme underperforms while modDG goes head to head with eRK4. modDG retains this same behavior even when we substitute simpler expression for the $\delta$ function.}
$\\[1ex]$

\end{minipage} \begin{minipage}[c]{0.01\textwidth}   \end{minipage} \begin{minipage}[c]{0.495\textwidth}

\includegraphics[width=6cm, height=5cm]{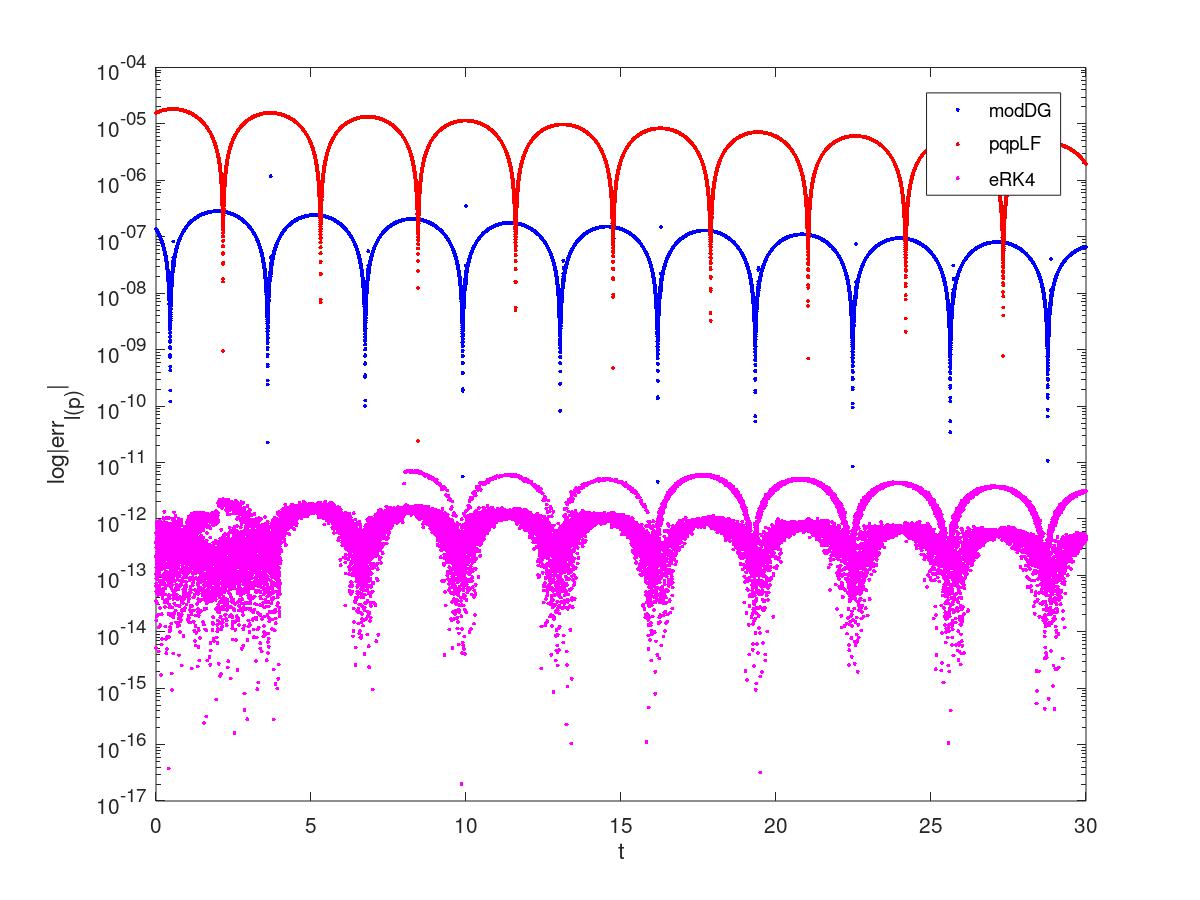}

\footnotesize{Figure 4.3.: Local error of $p$. For $q$ our scheme and SV are both committing errors of order magnitude $10^{-7}$.}
$\\[1ex]$

\end{minipage}

Minor investigation is in order to measure qualitative features of this new scheme with respect to others, e.g.: initial energy preservation (clearly better, as we saw) and energy decrementation rate.

We introduce the quantity
\be
R=\frac{E_{i+1}}{E_i}
\ee
describing the energy loss ratio of the system.

\begin{minipage}[c]{0.495\textwidth}

\includegraphics[width=6cm, height=5cm]{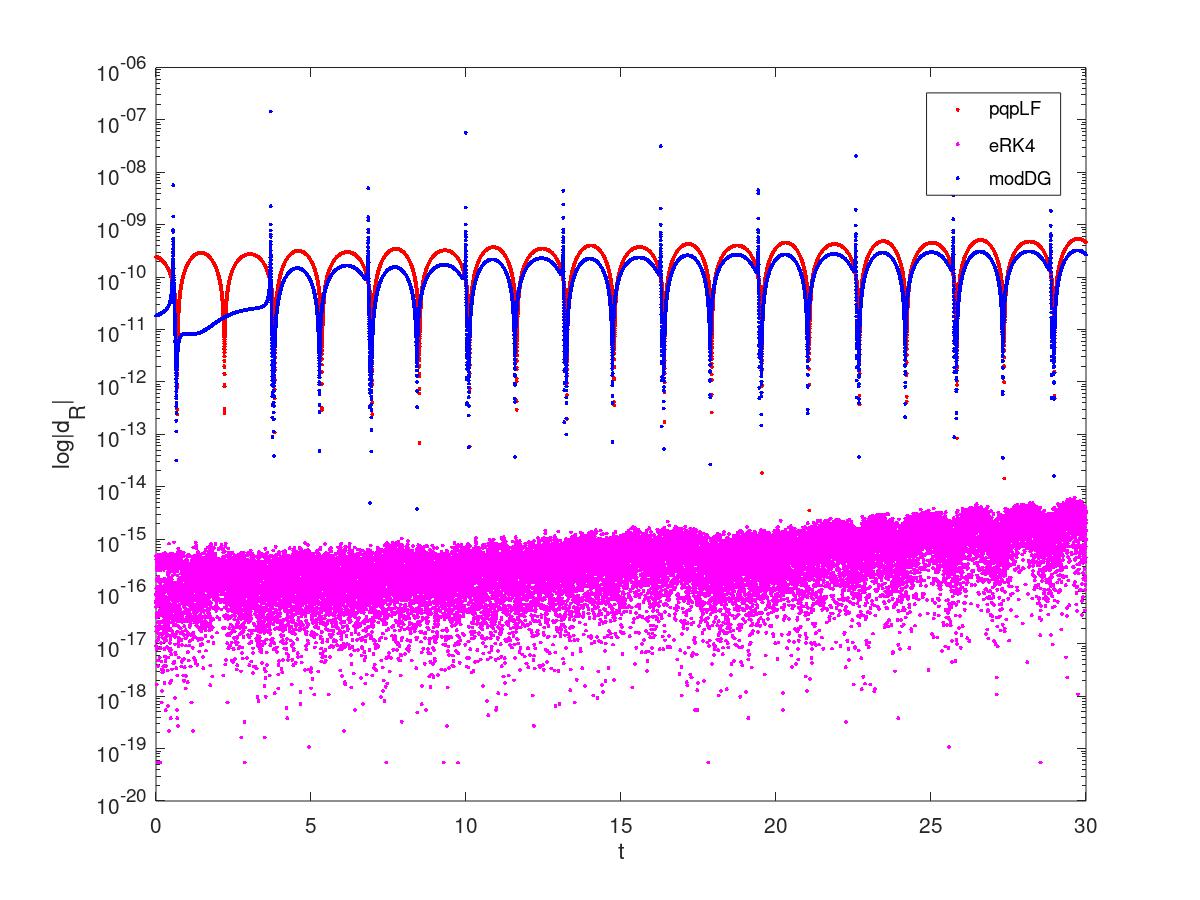}

\footnotesize{Figure 4.4.: $d_R$ treated by both fourth order schemes and pqpLF.}
$\\[1ex]$

\end{minipage} \begin{minipage}[c]{0.01\textwidth}   \end{minipage} \begin{minipage}[c]{0.495\textwidth}

This ends our investigation for numerical proof of non-potential Hamiltonian mechanics being effective. We have confirmed new differentiation rules being in perfect agreement with symplectic schemes like pqpLF, we witnessed new discrete gradient being proper method for dealing with dissipative systems, but more work have to be done in order to adjust accurracy to more delicate tasks.
$\\[1ex]$

\end{minipage}

\section{Conclusions and plans}

Inside the teritory of GNI we seemed always to work in favour of conservative systems, and results we obtained in this paper shed new light on the matter: dissipative systems may also be treated along the lines of gradient methods and symplectic schemes. This is partialy due to introduction of effectively conserved quantity $K$.

In the future we must take into account that the framework presented here suffers from many formal issues: The main flaw is that the Poisson bracket is not the entity that works well with reservoirs, the Jacobi identity breaks down and canonical Poisson bracket $\{q,p\}$ should not be equal to one, if Jacobi identity is to be saved. Additionally, for robust use of phase space techniques, uniqueness problem for phase trajectories have to be reconciled, as mentioned in the text.

Numerically, we saw the problem with rising accuracy.

Nevertheless, we can be happy with what was achieved: systematic treatment was proposed and it did not fail to accomplish given objectives; simulations were performed and their results will be published elsewhere for the Duffing oscillator, Van der Pol oscillator (both pure and modified) by both discrete gradient and symplectic counterparts, enabling us to develop further on the subject of classical energy reservoirs and pointing in the direction of new interesting numerical concepts.

\end{document}